# General Lower Bounds based on Computer Generated Higher Order Expansions

**Martijn A.R. Leisink and Hilbert J. Kappen**
University of Nijmegen, Department of Biophysics
Geert Grooteplein 21, 6525EZ Nijmegen, The Netherlands
{martijn,bert}@mbfys.kun.nl

## Abstract

In this article we show the rough outline of a computer algorithm to generate lower bounds on the exponential function of (in principle) arbitrary precision. We implemented this to generate all necessary analytic terms for the Boltzmann machine partition function thus leading to lower bounds of any order. It turns out that the extra variational parameters can be optimized analytically. We show that bounds upto nineth order are still reasonably calculable in practical situations. The generated terms can also be used as extra correction terms (beyond TAP) in mean field expansions.

## 1 INTRODUCTION

Mean field based expansions among which TAP (second order) in (Thouless et al., 1977) and (Plefka, 1981) and third order in (Leisink and Kappen, 2001) have in common that there is a need to compute higher order expansion terms. For the first few orders it is already an elaborate process to find the analytic expressions, but for higher orders the help of a computer algorithm is indispensible. With the computer speed nowaday, it is possible to generate these terms and, more importantly, to compute their contribution. For practical problems this often leads to a better approximation. The goal of this article is to develop an algorithm that computes a lower bound on the partition function of a Boltzmann machine (Ackley et al., 1985) of arbitrary precision, only limited by the patience of the user.

In (Leisink and Kappen, 2001) it was shown how one can find all polynomials of any (odd) order, which are lower bounds on the exponential function. In this article we make this rather theoretical idea applicable to real world problems by implementing an algorithm that generates lower bounds of arbitrary order on the Boltzmann machine partition function. These can be used either to compute more accurate approximations for means and correlations or, in combination with upper bounds as in (Jaakkola and Jordan, 1996), to restrict these statistics to certain regions (Leisink and Kappen, 2002).

In section 2 we briefly recall the algorithm to obtain the coefficients for the bounding polynomials. It turns out that for each two orders, there is one variational parameter which can be chosen freely. In section 3 we show that the optimal values for these variational parameters can easily be found. A time consuming optimization algorithm is not needed. Moreover, we will prove in that section that increasing the order never makes the bound less tight.

After these general theoretical considerations, the framework is applied to the Boltzmann machine in section 4. This step is far from trivial. We need to find all possible ways that two neurons can couple in the analytic expansion. Although this leads to an awful lot of expansion terms, we will show that the computational complexity is still reasonable. For instance, the time to compute a nineth order bound scales with network size as $\mathcal{O}(N^4)$. In section 5, we show the results of several numerical simulations and finally, in section 6, we discuss the applicability of the general bound to graphical models other than the Boltzmann machine.

## 2 THE CLASS OF LOWER BOUNDING POLYNOMIALS

In (Leisink and Kappen, 2001) it is shown how we can use a known bound on a function to obtain higher order bounds. The procedure is as follows: Given that $F_0(x) \geq B_0(x)$ we can create two primitive functions $F_1(x) = \int F_0(x)$ and $B_1(x) = \int B_0(x)$ such that $F_1(\mu) = B_1(\mu)$ for some $\mu$. If we apply this procedure again, thus constructing primitive functions $F_2(x) = \int F_1(x)$ and $B_2(x) = \int B_1(x)$ such that they are equal at the same point $x = \mu$, one can prove that

$F_2(x) \geq B_2(x)$ for all $x$. Specifically, for the exponential function, this yields

$$F_0(x) = e^x \geq 0 = B_0(x)$$
$$F_1(x) = e^x \text{ and } e^\mu = B_1(x)$$
$$\forall_{\mu,x} \, F_2(x) = e^x \geq e^\mu(1-\mu) + e^\mu x = B_2(x) \quad (1)$$

$B_2(x)$ is the well known tangential bound with the touching point at $x = \mu$. This series can be continued to obtain higher order bounds, which all can be written as

$$e^x \geq B_K(x) = \sum_{n=0}^{K-1} A_{K;n} x^n \quad (K \text{ even}) \quad (2)$$

where $A_{K;n}$ are coefficients of the polynomial in $x$.

Instead of writing down the higher order bounds explicitly, we can define the coefficients by the following recursive relation for $k = 0, 2, 4 \ldots K-2$:

$$\begin{cases} A_{k+2;n+2} = \dfrac{A_{k;n}}{(n+2)(n+1)} \quad \forall_{n \geq 0} \\ A_{k+2;1} = e^{\mu_k} - \sum_n \dfrac{A_{k;n}}{n+1} \mu_k^{n+1} \\ A_{k+2;0} = e^{\mu_k}(1 - \mu_k) + \sum_n \dfrac{A_{k;n}}{n+2} \mu_k^{n+2} \end{cases} \quad (3)$$

Note that coefficients $A_{k;n}$ and variational parameters $\mu_k$ only exist for $k$ is even since only odd order polynomials can bound the exponential function tightly. Starting with $\forall_n A_{0;n} = 0$ we can find all possible polynomial bounds by evaluating the recursive relation for $k = 0$, $k = 2$, etc. upto $k = K-2$. The above recursive relation is valid for all non-negative $n$, but as a consequence of the starting conditions, $A_{k;n} = 0$ for all $n \geq k$. Thus, finally, we have the coefficients $A_{K;n}$ for $n = 0 \ldots K-1$ and the variational parameters $\mu_k$ for $k = 0, 2, 4 \ldots K-2$, which together define the $K-1$-st order bound, $B_K(x)$.

The third order bound, for example, is given by

$$B_4(x) = \left[ e^{\mu_2}(1 - \mu_2) + \frac{1}{6} e^{\mu_0}(3 - 2\mu_0)\mu_0^2 \right] +$$
$$\left[ e^{\mu_2} - \frac{1}{2} e^{\mu_0}(2 - \mu_0)\mu_0 \right] x + \quad (4)$$
$$\left[ \frac{1}{2} e^{\mu_0}(1 - \mu_0) \right] x^2 + \left[ \frac{1}{6} e^{\mu_0} \right] x^3$$

where the square brackets are the coefficients $A_{4;0}$, $A_{4;1}$, $A_{4;2}$ and $A_{4;3}$, respectively. These coefficients are functions of the variational parameters $\mu_0$ and $\mu_2$, which can take any value without violating the bounding property.

## 3 OPTIMIZED BOUNDS FOR GRAPHICAL MODELS

For many graphical models, the log-probability of finding it in a state $\vec{s}$, is proportional to some energy function, thus $p(\vec{s}) \propto \exp(H(\vec{s}))$. A common problem is computing the normalizing constant of the distribution $p(\vec{s})$, since this requires the summation over exponentially many terms. Fortunately, we can use the bound derived above to lower bound the normalizing function $Z$:

$$Z = \sum_{\text{all } \vec{s}} \exp(H(\vec{s})) \geq \sum_{\text{all } \vec{s}} e^{\tilde{H}(\vec{s})} B_K \left( H(\vec{s}) - \tilde{H}(\vec{s}) \right)$$
$$= \tilde{Z} \sum_{n=0}^{K-1} A_{K;n} \langle \Delta H^n \rangle \quad (5)$$

where

$$\tilde{Z} = \sum_{\text{all } \vec{s}} e^{\tilde{H}(\vec{s})} \quad (6)$$

and $\langle \cdot \rangle$ denotes an average over the probability distribution with energy function $\tilde{H}(\vec{s})$. $\Delta H$ is an abbreviation for $H(\vec{s}) - \tilde{H}(\vec{s})$. The bound is valid for any $\tilde{H}(\vec{s})$, but obviously there is the constraint that the right hand side should be tractable to compute. This equation is the most general form for bounding the partition function with an odd order polynomial. Note that in general the variational parameters itself are allowed to be functions of $\vec{s}$. In this article, however, we will assume them to be constants, such that in equation 5 the coefficients $A_{K;n}$ can be taken out of the average.

To find the tightest bound, we set all variational parameters $\mu_i$ such that the bound is maximized. Hence the optimal $\mu_i$ satisfy

$$\frac{\partial}{\partial \mu_i} \left( \tilde{Z} \cdot \langle B_K(\Delta H) \rangle \right) = \tilde{Z} \cdot \left\langle \frac{\partial B_K}{\partial \mu_i} \right\rangle = 0 \quad (7)$$

It might be unexpected, but we can directly find the solution of equation 7. There is no need to apply any kind of maximization algorithm, which we will explain now.

Instead of taking the derivative with respect to $\mu_i$ of the explicit expression of $B_K$, we can perform this operation on the recursive relation defined in equation 3, which yields:

$$\begin{cases} A'_{k+2;n+2} = \dfrac{A'_{k;n}}{(n+2)(n+1)} \quad \forall_{n \geq 0} \\ A'_{k+2;1} = E_i \delta_{ik} - \sum_n \dfrac{A'_{k;n}}{n+1} \mu_k^{n+1} \\ A'_{k+2;0} = -\mu_i E_i \delta_{ik} + \sum_n \dfrac{A'_{k;n}}{n+2} \mu_k^{n+2} \end{cases} \quad (8)$$

where the prime denotes differentiation with respect to $\mu_i$ and

$$E_i = e^{\mu_i} - \sum_n A_{i;n}\mu_i^n \qquad (9)$$

Starting with $\forall_n A'_{0;n} = 0$ and evaluating this recursive relation for $k = 0$, $k = 2$, etc. upto $k = K-2$ defines the differentiated bound

$$\left\langle \frac{\partial B_K}{\partial \mu_i} \right\rangle = \sum_{n=0}^{K-1} A'_{K;n} \langle \Delta H^n \rangle \qquad (10)$$

One important property directly follows from equation 8. When evaluating it starting with $A'_{0;n} = 0$ and thus generating the coefficients of the differentiated bound, the first non-zero coefficient enters the equation precisely at the moment that $k = i$. As a consequence, we can write down explicit expressions for $A'_{i+2;n}$:

$$\begin{cases} A'_{i+2;n+2} &= 0 \\ A'_{i+2;1} &= E_i \\ A'_{i+2;0} &= -\mu_i E_i \end{cases} \qquad (11)$$

Secondly, we notice that the only dependency between $A'_{i+2;n}$ and $\mu_j$ with $j < i$ is through the function $E_i$.

At this point we define $X_i(x)$ to be the polynomial for which the coefficients are given by evaluating recursive relation 8 from $k = i+2$ upto $k = K-2$ starting with $\forall_{n \neq 1} A'_{i+2;n} = 0$ and $A'_{i+2;1} = 1$. $Y_i(x)$ is defined similarly, but starting with $\forall_{n \neq 0} A'_{i+2;n} = 0$ and $A'_{i+2;0} = -1$. Thus we may write

$$\left\langle \frac{\partial B_K}{\partial \mu_i} \right\rangle = E_i \left\langle X_i(\Delta H) + \mu_i Y_i(\Delta H) \right\rangle = 0 \qquad (12)$$

Now, regardless the exact value of $E_i$ in this equation the optimal value for $\mu_i$ is given by

$$\mu_i^{\text{opt}} = -\frac{\langle X_i(\Delta H) \rangle}{\langle Y_i(\Delta H) \rangle} \qquad (13)$$

Since $X_i$ and $Y_i$ do not depend on $\mu_j$ with $j \leq i$, the optimal value of $\mu_i$ only depends on $\mu_j$ with $j > i$. Therefore we can start computing the last variational parameter, $\mu_{K-2}^{\text{opt}}$, use that to compute $\mu_{K-4}^{\text{opt}}$, etc. We find, for example, that

$$\mu_{K-2}^{\text{opt}} = \langle \Delta H \rangle \qquad (14)$$

$$\mu_{K-4}^{\text{opt}} = \frac{\frac{1}{6}\langle \Delta H^3 \rangle - \frac{1}{2}\langle \Delta H \rangle \mu_{K-2}^{\text{opt}\,2} + \frac{1}{3}\mu_{K-2}^{\text{opt}\,3}}{\frac{1}{2}\langle \Delta H^2 \rangle - \langle \Delta H \rangle \mu_{K-2}^{\text{opt}} + \frac{1}{2}\mu_{K-2}^{\text{opt}\,2}} \qquad (15)$$

$\vdots$

All important quantities, such as $\mu_i^{\text{opt}}$ and $\langle B_K(\Delta H) \rangle$, can easily be computed by a computer program using the recursive relations shown in section 2 and 3. Therefore, there is no need to write down the full analytic expressions for any order, while still being able to compute their value.

The attentive reader might have thought about choosing $\mu_i$ such that $E_i$ in equation 12 becomes zero, which gives rise to another solution for $\mu_i^{\text{opt}}$. This solution, however, can never correspond to a maximum. Looking at equation 9, we see that $E_i$ is in fact the difference between $e^{\mu_i}$ and the bound $B_i(x)$ evaluated at $x = \mu_i$. Therefore it is obvious that $E_i \geq 0$. This implies that the right hand side of equation 12 will not change its sign from plus to minus when it passes through $E_i = 0$. Since this behaviour is a requirement for a point to be a maximum, this solution should not be taken[1].

In the appendix, it will be shown that the polynomial $Y_i(\Delta H)$ is always negative. Therefore, the derivative at the point $\mu_i^{\text{opt}}$ as given in equation 13 does change the sign from positive to negative and thus corresponds to a maximum.

It is worth mentioning that when all $\mu_i$ are set to zero, the bound $B_K(x)$ coincides with the Taylor expansion of the exponential function around zero upto the $K{-}1$-st order. Therefore the lower bound expansion has an infinite radius of convergence. This is in contrast with the Plefka expansion in (Plefka, 1981), where it is proven that such an expansion (e.g. TAP) suffers from a finite radius of convergence. Therefore computing extra correction terms for the Plefka expansion outside this radius of convergence does in general not make the results more accurate. The lower bound expansion, on the other hand, can approximate the partition function with any desired accuracy, although this is, of course, strongly limited by the available computer time. One word of caution: One has to ensure that the distribution $\tilde{H}(\vec{s})$ is such, that none of the exponentially many terms in equation 5 has a contribution which is smashed to zero, because of a nearly zero $\exp(\tilde{H}(\vec{s}))$ term (i.e. $\tilde{H}(\vec{s})$ should represent a distribution which is 'flat enough'). At the end of section 5 we will briefly show some results that show this effect.

To conclude this section, we remark that whenever we have an optimized bound $\langle B_K(\Delta H) \rangle$, we can construct a next order bound $\langle B_{K+2}(\Delta H) \rangle$ which is identical. For this purpose, we set $\mu_0 = -\infty$ for the second bound and set its $\mu_{k+2}$ to the $\mu_k^{\text{opt}}$ of the first one. It is immediately clear from the definitions in section 2 that these bounds are identicle. The bound $\langle B_{K+2}(\Delta H) \rangle$, however, is not necessarily in a maximum and can be optimized further. Therefore, $\langle B_{K+2}(\Delta H) \rangle$ is at least as tight as $\langle B_K(\Delta H) \rangle$.

---

[1]In fact, the solution for $\mu_i$ where $E_i = 0$ corresponds to a point, where the bound has a shape similar to the function $y = x^3$ at zero.

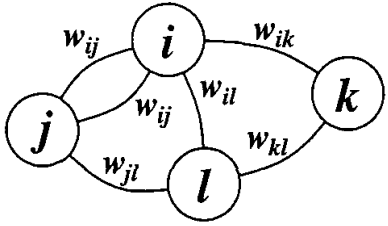

Figure 1: A visualisation of a sixth order partition given by $w_{ij}w_{ij}w_{ik}w_{il}w_{jl}w_{kl}$.

## 4 BOLTZMANN MACHINES

Boltzmann machines fit into the general framework, in which their energy function is defined by $H(\vec{s}) = \frac{1}{2}\sum_{ij} w_{ij} s_i s_j + \sum_i \theta_i s_i$, where $w_{ij}$ are symmetric weights and $\theta_i$ thresholds on binary valued neurons. The distribution defined by $\tilde{H}(\vec{s})$ is often chosen to be factorized[2], thus $\tilde{H}(\vec{s}) = \sum_i h_i s_i + constant$. These definitions allows us to write the terms in equation 5 as

$$\langle \Delta H^n \rangle = \left\langle \left( \frac{1}{2}\sum_{ij} w_{ij}(s_i s_j - m_i m_j) + \sum_i (\theta_i - h_i)(s_i - m_i) \right)^n \right\rangle \quad (16)$$

where the constant is chosen such that $\langle \Delta H \rangle = 0$ and $m_i = \langle s_i \rangle = \tanh h_i$.

To simplify the computations dramatically, we require that the $h_i$'s obey the mean field equations:

$$\forall_i \; h_i = \theta_i + \sum_j w_{ij} m_j \quad (17)$$

Given this property of $h_i$ and the fact that the weights are symmetric we can rewrite equation 16 as

$$\langle \Delta H^n \rangle = \left\langle \left( \frac{1}{2}\sum_{ij} w_{ij}(s_i - m_i)(s_j - m_j) \right)^n \right\rangle \quad (18)$$

Although this expression does not look very hard to compute, the opposite is true. This is due to the non-trivial coupling of neurons with equal indices. Whenever a pair $s_i s_i$ occurs, this evaluates to the constant one instead of $m_i^2$. Therefore we need to find all possible ways the indices can couple, each of them called a 'partition'. For example, all second order partitions are $w_{ij}w_{kl}$, $w_{ij}w_{ik}$ and $w_{ij}w_{ij}$.

---
[2]Other choices, which preserve the tractability of the bound, are possible.

It turns out to be useful to make a correspondance between a partition and a graph. This can be done by drawing as many nodes as we have independent indices and draw links between them if there is a weight having these two indices. For instance, the partition $w_{ij}w_{ij}w_{ik}w_{il}w_{jl}w_{kl}$, which can occur for $n = 6$, can be visualised as in figure 1. The contribution of this partition to equation 18 can easily be computed using the graph. Firstly, each node is assigned the vector $M_c = \langle (s_x - m_x)^c \rangle$, where $c$ is the number of connections to that particular node and $x$ stands for the corresponding index that node is representing. This term can be seen as the $c$-th moment of the factorized distribution (see also table 1). After that, we can apply a kind of junction tree algorithm (or, similarly, a variable elimination scheme) to compute the contribution of this partition. See (Lauritzen and Spiegelhalter, 1988) for a detailed description of the junction tree algorithm. In the example we start by summing out the index $k$, then $l$, and finally $i$ and $j$ thus leading to the value of the average over this partition.

The only problem that frustrates the independently summing out of indices, is the requirement that no pair of indices should take equal values (since that contribution is captured in another partition). We can, however, let all indices run freely as long as we correct for this elsewhere. This can be done in the following way: To all nodes with a single connection, we assign the vector $M_1$. All nodes with two connections are assigned the vector $M_2 - M_1^2$, where $M_1^2$ is the correction for the fact that the two indices in more refined partitions were allowed to be equal. These corrected moments are written with a prime as $M_c'$. See table 1 for more examples.

Thus the full expression to be computed for the partition in figure 1 reads as

$$\sum_{ij} \left( w_{ij}^2 \left(-2 + 8m_i^2 - 6m_i^4\right) \left(-2m_j + m_j^3\right) \right.$$
$$\sum_l \left( w_{il} w_{jl} \left(-2m_l + 2m_l^3\right) \right.$$
$$\left. \left. \sum_k \left( w_{ik} w_{kl} \left(1 - m_k^2\right) \right) \right) \right) \quad (19)$$

where the terms are grouped such that it can be computed the most efficiently. In this case that is a computational complexity proportional to $N^3$ (with $N$ the number of neurons of the Boltzmann machine). Obviously, three is also the maximum clique size of the example partition graph in figure 1.

Since the maximum number of couplings in any partition graph is equal to the expansion order of the bound, the computational complexity of an expansion upto order $n$ scales as $N^{\pi(n)}$, where $\pi(n)$ is the size of the largest clique one can build with $n$ couplings. This is roughly equal to $\sqrt{2n}$. This implies that a nineth order bound, for instance, scales with $N^4$, since one needs at least ten couplings to construct a clique of size five.

Table 1: Moments and corrected moments. The latter can be computed by taking the moment $M_c$ and substracting all possible combinations of corrected moments above that line which subscripts add up again to $c$. Obviously, the coefficient in front of each of the correction terms is calculated as the number of possible partitions of $c$ into the specified subsets. Since $M_1 = 0$ it is not written down.

| Moments | Corrected moments |
|---|---|
| $M_1 = \langle s_i - m_i \rangle = 0$ | $M_1' = 0$ |
| $M_2 = \langle (s_i - m_i)^2 \rangle = 1 - m_i^2$ | $M_2' = M_2 = 1 - m_i^2$ |
| $M_3 = \langle (s_i - m_i)^3 \rangle = -2m_i + 2m_i^3$ | $M_3' = M_3 = -2m_i + 2m_i^3$ |
| $M_4 = \langle (s_i - m_i)^4 \rangle = 1 + 2m_i^2 - 3m_i^4$ | $M_4' = M_4 - 3M_2'^2 = -2 + 8m_i^2 - 6m_i^4$ |
| $M_5 = \langle (s_i - m_i)^5 \rangle = -4m_i + 4m_i^5$ | $M_5' = M_5 - 10M_2'M_3' = 16m_i - 40m_i^3 + 24m_i^5$ |
| $M_6 = \langle (s_i - m_i)^6 \rangle$ | $M_6' = M_6 - 15M_2'M_4' - 10M_3'^2 - 15M_2'^3$ |
| $= 1 + 9m_i^2 - 5m_i^4 - 5m_i^6$ | $= 16 - 136m_i^2 + 240m_i^4 - 120m_i^6$ |
| $M_7 = \langle (s_i - m_i)^7 \rangle$ | $M_7' = M_7 - 21M_2'M_5' - 35M_3'M_4' - 105M_2'^2 M_3'$ |
| $= -6m_i - 14m_i^3 + 14m_i^5 + 6m_i^7$ | $= -272m_i + 1232m_i^3 - 1680m_i^5 + 720m_i^7$ |

The final step in computing expression 18 is finding all possible partitions together with how many times each of them occurs. It is not obvious how to do this search efficiently (comparing graphs is NP-hard), but it is possible to develop quite fast algorithms for this purpose. Fortunately, these results are problem independent. Thus one can compute them once upto some order and store them forever. The final number of distinct partitions is reasonable, although it scales at least exponentially with $n$. In table 2 all partition graphs are explicitly shown for the first five orders of the bound. In table 3 the number of times each partition graph occurs is shown for the first five orders. On http://www.mbfys.kun.nl/~martijn one can find all partition graphs together with the program that generated them.

## 5 NUMERICAL RESULTS

To assess the quality of the any order bounds, we generated 9,000 fully connected Boltzmann machines with $N = 14$ neurons. The $\theta_i$ were drawn from a Gaussian with zero mean and standard deviation 0.2 and $w_{ij}$ drawn with standard deviation $\sigma_w/\sqrt{N}$. This is known as the SK-model as in (Sherrington and Kirkpatrick, 1975). For all networks we computed the exact partition function (which is still tractable in this regime) and the lower bound for several orders. Note that the 17th order bound is probably too time consuming to be computed for larger networks. Here it is only shown to illustrate the behaviour of the really high order bounds. Nineth order, however, is perfectly doable (see table 2).

Table 2: Computational complexity of the expansions. The actual complexity is slightly higher, since only the partitions of maximum clique size are taken into account (the leading term). The sum between brackets shows how many of all the partitions had clique size two, three, etc. Partitions in which a node with a single connection occurs, are not counted, since $M_1' = 0$ anyway. Note that although the number of partitions scales rather badly with the order, the scaling with network size is reasonably small.

| Order | # partitions | $\pi(Order)$ | Complexity |
|---|---|---|---|
| 2 | 1 (1) | 2 | $\sim \mathcal{O}(N^2)$ |
| 3 | 2 (1+1) | 3 | $\sim \mathcal{O}(N^3)$ |
| 4 | 5 (3+2) | 3 | $\sim \mathcal{O}(2 \cdot N^3)$ |
| 5 | 11 (4+7) | 3 | $\sim \mathcal{O}(7 \cdot N^3)$ |
| 6 | 34 (11+22+1) | 4 | $\sim \mathcal{O}(N^4)$ |
| 7 | 87 (18+67+2) | 4 | $\sim \mathcal{O}(2 \cdot N^4)$ |
| 8 | 279 (45+221+13) | 4 | $\sim \mathcal{O}(13 \cdot N^4)$ |
| 9 | 897 (91+744+62) | 4 | $\sim \mathcal{O}(62 \cdot N^4)$ |

Table 3: All graphs with a non-zero contribution belonging to second, third, fourth and fifth order expansion terms. The numbers indicate how many times they occur.

| |
|---|
| $2 \cdot$ ⬤—⬤ |
| $4 \cdot$ ◿ $+ 8 \cdot$ △ |
| $8 \cdot$ ◿ $+ 96 \cdot$ △ $+ 48 \cdot$ ⋈ $+ 12 \cdot$ ⬤⬤ $+ 48 \cdot$ ◻ |
| $16 \cdot$ ◿ $+ 320 \cdot$ △ $+ 480 \cdot$ △ $+ 320 \cdot$ ⋈ $+$ $80 \cdot$ ⬤⬤ $+ 480 \cdot$ ⬤⬤ $+ 960 \cdot$ ◻ $+ 960 \cdot$ ◻ $+$ $960 \cdot$ ◻ $+ 160 \cdot$ ⬤⬤ $+ 384 \cdot$ ⬤ |

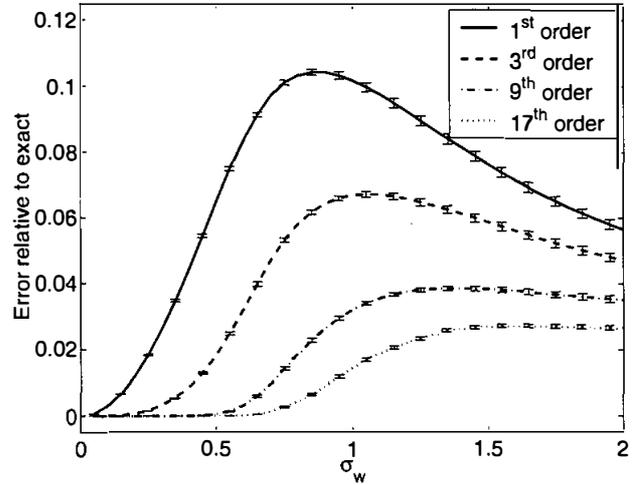

Figure 2: The relative error $\mathcal{E}$ of the bound compared to the true log partition function for various sizes of the weights. We computed the mean error for a set of 150 networks. Plotted are the mean and standard deviations for 60 of such means.

The relative error, $\mathcal{E} = 1 - \log B_K / \log Z$, of the optimized bound compared to the exact partition function is shown in figure 2. The error bars are standard deviations of the means, thus indicating that the mean curves are significantly different. This gives no information about whether the tightness of the bound would increase with expansion order for only one instance of a network. This is, however, a theoretical fact (see section 3).

Besides these bounding errors, it is reasonable to suspect that a better approximation of the partition function also leads to more accurate means and correlations, which are derived quantities. For the same 9,000 networks, we computed the exact correlation between the first two neurons and the approximated ones using the bounds as approximations for the partition function. It is clear from figure 3 that indeed the correlations are more accurately computed by using higher order bounds. Unfortunately, the so obtained values are neither upper nor lower bounds on the correlations. One could, however, combine the improved lower bounds with already existing upper bounds for Boltzmann machines to find definite regions in which the means and correlations must lie.

In case of the bounding error (figure 2), it is obvious that all errors tend to zero for small weight sizes. For large weights, however, the errors become closer and closer to each other as well. This can be understood since for very large weights, there is usually one eigenvector of the weight matrix, that overwhelms the others, thus leading to a single pair of opposite states, which has non-zero probability. This can perfectly be catched by a factorized model. The drawback of this, however, is that the mean field solution for $h_i$ starts heading towards plus and minus infinity much earlier, which surpresses the contribution of other terms than those two states (see also the remarks at the end of section 3). This could be the optimal strategy in case of a first order (or mean field) bound, but this does not necessarily hold for higher orders. This might be solved by trying to find other solutions of the mean field equations than the standard ones[3]. For instance, $h_i = 0$ is always a solution if the thresholds are zero, although this solution is usually not found when solving equation 17 iteratively (this point is only an attractor for weights with $\sigma_w$ less than about a half).

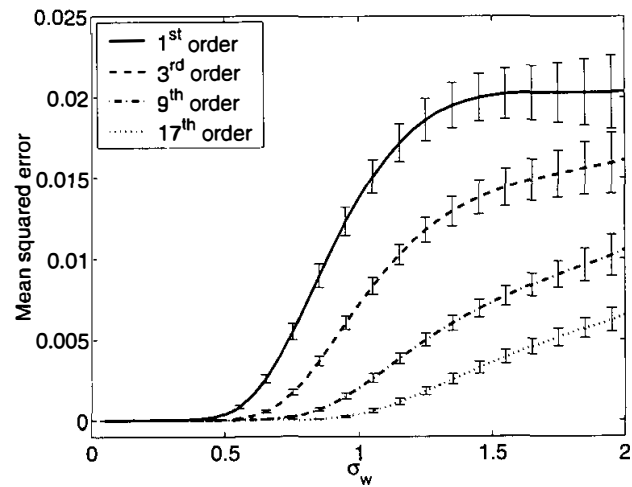

Figure 3: The mean squared error of the correlations following from the bound and the exact ones. We computed the mean error for a set of 150 networks. Plotted are the mean and standard deviations for 60 of such means.

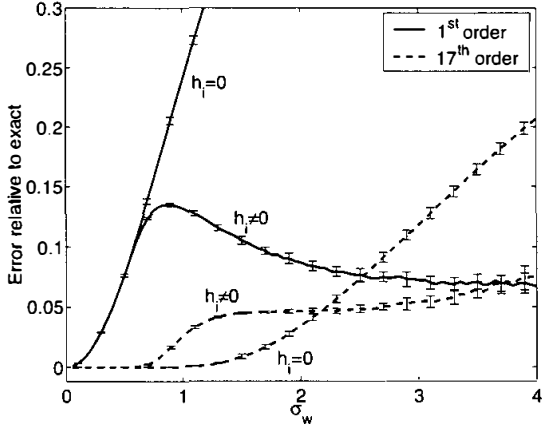

Figure 4: The relative error of the bound compared to the true log partition function for various sizes of the weights. All thresholds were zero. The error for the first and 17th order bound are shown given the ordinary mean field solution $h_i \neq 0$ and the (possibly non-stable) solution $h_i = 0$, which also obeys equation 17. We computed the mean error for a set of 100 networks. Plotted are the mean and standard deviations for 30 of such means.

In figure 4 the effect of taking another mean field solution is shown for a Boltzmann machine with zero thresholds. Obviously, the first order bound becomes worse by taking $h_i = 0$ as a solution, since the standard solution of equation 17 follows directly from optimizing this bound. The 17th order bound, however, is much tighter in a certain region for this (unstable) mean field solution. A similar behaviour can be seen for the other higher order bounds (not shown here). A general procedure to find such solutions is still a topic of research.

---

[3]We could, of course, take any $h_i$, but a solution to the mean field equations has the advantage that is simplifies equation 16 enormously.

## 6 APPLICABILITY TO OTHER GRAPHICAL MODELS

Although section 4 specifically showed the applicability of the any order bounds to the Boltzmann machine, the algorithm outlined in section 2 and 3 can be used for a larger class of graphical models. Generally speaking, it can be applied to any model for which we can find a distribution with energy function $\tilde{H}(\vec{s})$ such that the average $\langle \Delta H^n \rangle$ under that distribution (see equation 5) is tractable to compute.

Markov networks, for instance, where the energy function is defined by the sum over a number of potentials over groups of neurons

$$H(\vec{s}) = \sum_{p=1}^{P} \psi_p(S_p) \quad (S_p \text{ is a subset of nodes}) \quad (20)$$

belongs to this class of networks. In this case, choosing $\tilde{H} = \sum_i h_i s_i + constant = \sum_i \tilde{\psi}_i$ (as for the Boltzmann machine) will make $\langle \Delta H^n \rangle$ tractable to compute as long as $n$ is not too high. At a first glance, the worst case seems to be that the number of terms one has to compute is about $P^n$, where $P$ is the number of potentials in equation 20 and each term contains $n$ times as many nodes as the maximum of a potential. This is, however, too pessimistic thinking.

Take, for instance, the second order term, where we assume that the $\tilde{\psi}_i$ are incorporated in the $\psi_p$:

$$\langle \Delta H^2 \rangle = \sum_{pq} \langle \psi_p \psi_q \rangle \quad (21)$$

$$= \sum_{pq} (\langle \psi_p \psi_q \rangle - \langle \psi_p \rangle \langle \psi_q \rangle) + \left(\sum_p \langle \psi_p \rangle\right)^2$$

The last term vanishes, since the constant in $\tilde{H}$ was chosen such that $\langle \Delta H \rangle = 0$. It is immediately clear from equation 21, that only potential pairs that overlap (i.e. are dependent) have a non-zero contribution to equation 21, since the averages are taken over a factorized distribution. This number is usually much smaller than computing all $P^2$ potential pairs. Similar results can be obtained for $n > 2$. This means that for a lot of real world Markov networks the dependence between $P$ and the computational complexity will be close to linear (at least for small $n$). Note that computing a single average over a product of $n$ potentials can still be time consuming due to the maximum clique size of the product, although this is usually not as worse as $n$ times the original clique size.

To illustrate this, we investigate the computational complexity for a Markov network with a structure as shown in figure 5. First of all, it is clear that the number of products of potentials that do connect scales linearly with the number of nodes and thus with $P$. Moreover, upto fifth order the maximum clique size of any product of potentials remains two. For sixth order the maximum occurring clique size is still only three (a cycle of six nodes). This is in sharp contrast with the worst case scenario mentioned earlier. The only number that increases very rapidly is how many different connecting products of potentials are possible per node. This number, however, only depends on the order of the expansion and not (similar to the number of partitions in the Boltzmann machine case) on the number of neurons. This makes the more general case of Markov networks tractable as well in many practical situations. A more thorough investigation will be topic of future research.

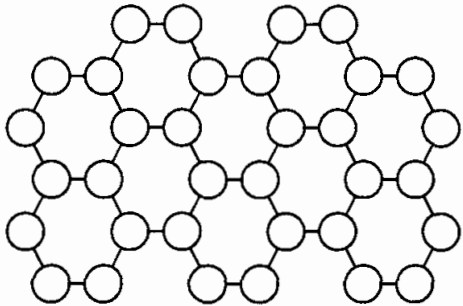

Figure 5: A network structure. Each line defines a potential of two nodes.

## APPENDIX

In this appendix we show that the polynomial $Y_i(\Delta H)$ in section 3 is always negative. When we start with the obvious bound $0 > -1$, we can construct a recursive relation very similar to equation 3, which yields

$$\begin{cases} A_{k+2;n+2} = \dfrac{A_{k;n}}{(n+2)(n+1)} \\ A_{k+2;1} = -\sum_n \dfrac{A_{k;n}}{n+1}\mu_k{}^{n+1} \\ A_{k+2;0} = \sum_n \dfrac{A_{k;n}}{n+2}\mu_k{}^{n+2} \end{cases} \quad (22)$$

In this case we start of course with $\forall_{n\neq 0} A_{0;n} = 0$ and $A_{0;0} = -1$. Thus given these starting conditions and a recursive relation as above, we find a lower bound on zero, or, in plain English, a negative number.

At this point we notice that the recursive relation (22) is identical to the one defined in equation 8 for any $k \neq i$. Moreover, $Y_i(\Delta H)$ was defined as the polynomial starting with $\forall_{n\neq 0} A'_{i;n} = 0$ and $A'_{i;0} = -1$ and applying the recursive relation (8) beginning with $k = i+2$. As we have seen in the previous paragraph, polynomials that are constructed in this way are in fact lower bounds on zero. Therefore, $Y_i(\Delta H) < 0$.

### Acknowledgements

This research is supported by the Technology Foundation STW, applied science devision of NWO and the technology programme of the Ministry of Economic Affairs.